\documentclass[12pt]{article}

\usepackage[top=50mm, bottom=50mm, left=50mm, right=50mm]{geometry}
\usepackage{lineno}
\usepackage{amssymb}
\usepackage{amsmath}
\usepackage{amsthm}
\usepackage{epsfig}
\usepackage{graphicx}
\usepackage{graphics}
\usepackage{float}
\usepackage{subfigure}
\usepackage{multirow}
\usepackage{color}
\usepackage{lineno}
\usepackage{fullpage}
\usepackage[normalem]{ulem} 
\usepackage{makeidx}
\usepackage{xspace}
\usepackage{wrapfig}
\usepackage{todonotes}
\usepackage{hyperref}
\usepackage{cleveref}
\usepackage{footnote}
\makeindex

\usepackage{accents}

\newtheorem{Definition}{Definition}[section]
\newtheorem{theorem}[Definition]{Theorem}

\theoremstyle{definition}
\newtheorem{Remark}[Definition]{Remark}

\definecolor{darkred}{rgb}{1, 0.1, 0.3}
\definecolor{darkblue}{rgb}{0.1, 0.1, 1}
\definecolor{darkgreen}{rgb}{0,0.6,0.5}

\newcommand {\mm}[1] {\ifmmode{#1}\else{\mbox{\(#1\)}}\fi}

\newcommand{\C}{\mathbb{C}} 
\newcommand{\T}{\mathbb{T}} 
\newcommand{\R}{\mathbb{R}} 
\newcommand{\Z}{\mathbb{Z}} 
\newcommand{\N}{\mathbb{N}} 
 
\newcommand{\hz}{\mathrm{HZ}}

\newcommand{\dd}{\mathrm{d}}

\begin{document}

\title{Symplectic capacities\\ of disc cotangent bundles of flat tori}
 
\author{Gabriele Benedetti\footnote{Department of Mathematics, Vrije Universiteit Amsterdam, De Boelelaan 1111, 1081 HV Amsterdam, g.benedetti@vu.nl}, Johanna Bimmermann\footnote{Fakult\"at f\"ur Mathematik, Ruhr-Universit\"at Bochum, Universit\"atsstra{\ss}e 150, D-44801 Bochum, Germany, Johanna.Bimmermann@rub.de, Kai.Zehmisch@rub.de}, and Kai Zehmisch\footnotemark[\value{footnote}]}

\maketitle

\begin{abstract}
\noindent
We show that on the unit disc cotangent bundle of flat Riemannian tori, all normalized capacities coincide with twice the systole. The same result holds for flat, reversible Finsler tori and normalized capacities that are greater than or equal to the Hofer--Zehnder capacity.
\end{abstract}

\section{Introduction}
Given $n\in\N$, symplectic capacities are invariants of symplectic manifolds $(M,\omega)$ of dimension $2n$ used to study symplectic embeddings between such manifolds. As introduced by Ekeland and Hofer in \cite{EH}, a (normalized) symplectic capacity is a function \[
c\colon\big\{\text{Symplectic manifolds of dimension $2n$}\big\}\to[0,\infty]
\]
satisfying the following properties:
\begin{description}
\item[Embedding] If $\varphi\colon(M_1,\omega_1)\to(M_2,\omega_2)$ is any symplectic embedding, then the inequality $c(M_1,\omega_1)\leq c(M_2,\omega_2)$ holds.
\item[Scaling] If $r>0$, then $c(M,r\omega)=rc(M,\omega)$.
\item[Normalization] If $\omega_0$ is the standard symplectic form on $\C^n$, then $c(B_r,\omega_0)=\pi r^2=c(Z_r,\omega_0)$, when $B_r$ is the open Euclidean ball of radius $r$ and $Z_r=D_r\times \C^{n-1}$ is the infinite symplectic cylinder having the open disc $D_r\subset\C$ of radius $r$ as base.
\end{description}
The first condition yields obstructions to the existence of symplectic embeddings. In view of this, the third condition encodes the celebrated Gromov's non-squeezing theorem: $B_{r_1}$ can be symplectically embedded in $Z_{r_2}$ if and only if $r_1\leq r_2$, see \cite{Gromov}.

For this reason, finding a function that satisfies the axioms of a capacity is a difficult problem. Nevertheless, a number of capacities have been discovered in the last forty years, see \cite{CHLS}. For instance, the Hofer--Zehnder capacity $c_{\hz}$ is defined using compactly supported Hamiltonian systems on $(M,\omega)$ and has a special dynamical significance: if $c_{\hz}(M,\omega)$ is finite then every proper Hamiltonian $H\colon M\to\R$ admits periodic orbits on almost all energy levels. See \cite[Ch.\ 4, Thm.\ 1]{HZ} and \cite{MS05} for a proof of this fact, and \cite{U12} for examples of closed manifolds having infinite Hofer--Zehnder capacity.\\

Despite these advances, it is a challenging task to compute capacities, even for simple classes of manifolds, or to show that two capacities coincide on these classes. For instance, the strong Viterbo conjecture asserts that all capacities coincide on the class of convex domains in $\C^n$ \cite[Section 14.9, Problem 53]{MS}. Another remarkable class of examples is given by canonical Finsler disc cotangent bundles $(D^*_FQ,\omega_\mathrm{can})$. Here,
\begin{itemize}
\item $Q$ is a closed manifold and
\[
D^*_FQ=\{(q,p)\in T^*Q\mid F^*(q,p)<1\},
\]
where $F\colon TQ\to\R$ is a smooth Finsler norm on $Q$ and $F^*\colon T^*Q\to\R$ denotes the dual norm
\[
F^*(q,p):=\sup_{\substack{v\in T_qQ\\ F(q,v)<1}}p\cdot v,\qquad \forall(q,p)\in T^*Q
\]
\item $\displaystyle\omega_\mathrm{can}=\sum_{i=1}^n\dd p_i\wedge\dd q_i$ is the canonical symplectic form.
\end{itemize}
Characteristics $(q,p)\colon\R\to \partial(D^*_FQ)$ project to Finsler geodesics $q\colon\R\to Q$, that is, extremals of the Finsler length
\[
\ell_F(q):=\int_a^bF(q(t),\dot q(t))\mathrm dt
\]
at fixed endpoints. We denote by $\mathrm{sys}(F)$ the Finsler length of the shortest non-contractible closed geodesic on $Q$. We call a Finsler norm reversible if $F(q,v)=F(q,-v)$ for all $(q,v)\in TQ$. The most fundamental examples of reversible Finsler norms are induced by Riemannian metrics $g$, that is, $F(q,v):=\sqrt{g_q(v,v)}$ for all $(q,v)\in TQ$.

Capacities of disc cotangent bundles are known only in some symmetric examples. The ball-capacity (also known as Gromov-width) of codisc bundles of $S^2$ and $\R\mathrm{P}^2$ with respect to the round metric was computed by Ferreira and Ramos \cite{FR22}. These results were also generalized to spheres of revolution in \cite{FRV23}. Furthermore, the Hofer--Zehnder capacity of disc cotangent bundles of $\C\mathrm{P}^n$ with respect to the Fubini--Study metric and of $\R\mathrm{P}^n$ with respect to the round metric was computed by the second author in \cite{B23.1}. Finally, in \cite{Jiang} Jiang computed the value of any normalized capacity for the codisc bundle of the $\ell^1$-norm on the $n$-torus. It is this last result that we aim to generalize in the present paper.

Indeed, our main theorem computes the capacity of disc cotangent bundles of flat, reversible Finsler norms on the torus $\T^n=\R^n/\Z^n$. We recall that a Finsler norm is flat if, in the global trivialization $T\T^n=\T^n\times\R^n$ induced by the quotient $\R^n\to\T^n$, there is $f\colon\R^n\to\R$ with $F(q,v)=f(v)$ for all $(q,v)\in \T^n\times \R^n$. In this case $D^*_F\T^n=\T^n\times K^*$, where $K^*=\{p\in(\R^n)^*\mid f^*(p)<1\}$.
\begin{theorem}\label{thm2}
If $F$ is a flat, reversible Finsler norm on $\T^n$ and $c$ is a (normalized) symplectic capacity, then
\[
c(D_F^*\T^n,\omega_\mathrm{can})\leq 2\mathrm{sys}(F).
\]
If, in addition,
\begin{enumerate}
    \item there is a Finsler norm $G$ induced by a flat Riemannian metric with $F\geq G$ and $\mathrm{sys}(F)=\mathrm{sys}(G)$, or
    \item $c\geq c_{\hz}$, where $c_{\hz}$ denotes the Hofer--Zehnder capacity, 
\end{enumerate}
then
\[
c(D_F^*\T^n,\omega_\mathrm{can})=2\mathrm{sys}(F).
\]
\end{theorem}
\begin{Remark}
In particular, $c(D_F^*\T^n,\omega_\mathrm{can})=2\mathrm{sys}(F)$ for all normalized capacities, if $F$ is induced by a Riemannian metric or if $F$ is an $\ell^p$-norm with $p\leq 2$. As recalled above, the case $p=1$ was proved in \cite{Jiang}.
\end{Remark}
Before ending this introduction, let us mention that it would be interesting to extend the computation of capacities to the larger class of twisted disc cotangent bundles $(D^*_FQ,\omega_\mathrm{can}+\pi^*\sigma)$, where $\pi\colon D^*_FQ\to Q$ is the footpoint projection, and $\sigma$ is a closed two-form on $Q$. In this case, characteristics on $\partial (D^*_FQ)$ correspond to magnetic geodesics which describe the motion of a particle in the stationary magnetic field $\sigma$. We refer the reader to \cite{B23,B23.1,B23.2} for computations of the Hofer--Zehnder capacity for twisted disc cotangent bundles of hermitian symmetric spaces, where the norm and the magnetic field are induced by the associated canonical K\"ahler structure.\\
\noindent \textbf{Acknowledgments.} G.B.~gratefully acknowledges support from the Simons Center for Geometry and Physics, Stony Brook University at which some of the research for this paper was performed during the program \textit{Mathematical Billiards: at the Crossroads of Dynamics, Geometry, Analysis, and Mathematical Physics}. G.B.~is partially supported by the DFG 
under Germany's Excellence Strategy EXC2181/1 - 390900948 (the Heidelberg STRUCTURES Excellence Cluster). J.B.~and K.Z.~are partially supported by the DFG under the Collaborative Research Center SFB/TRR 191 - 281071066 (Symplectic Structures in Geometry, Algebra and Dynamics). We are grateful to Alberto Abbondandolo, Hansjörg Geiges, Felix Schlenk and Stefan Suhr for inspiring conversations and their valuable suggestions on this work. We are indebted to the anonymous referee for their precious comments on the first draft of the manuscript.

\section{Proof of Theorem \ref{thm2}}
Let us identify $\T^n$ with the quotient of $\R^n$ by the standard lattice $\Z^n$. Geodesics of a flat Finsler norm $F$ are straight lines and therefore
\begin{equation}\label{sys}
s:=\mathrm{sys}(F)=\min\{f(v)\mid v\in\Z^n,\ v\neq0\},    
\end{equation}
where $F(q,v)=f(v)$ for all $(q,v)\in \T^n\times \R^n$.
We write $K=\{v\in\R^n\mid f(v)<1\}$ and $K^*=\{p\in(\R^n)^*\mid f^*(p)<1\}$, where $f^*$ is the dual Finsler norm, which, by reversibility, has the formula $f^*(p)=\sup_{v\in K}|p\cdot v|$.
\subsection{Upper bound}
Let $e_1,\ldots,e_n$ be the standard basis of $\Z^n$. Let $u=(k_1,\ldots,k_n)\in \Z^n$ be such that \[
f(u)=\min\{f(v)\mid v\in\Z^n,\ v\neq0\}=s.\]
In particular, the numbers $k_1,\ldots,k_n$ are coprime and, therefore, there exists a linear isomorphism $A\colon\R^n\to\R^n$ represented by a matrix in $SL(n,\Z)$ such that $A(e_1)=u$. This isomorphism induces a diffeomorphism $[A]\colon\T^n\to\T^n$ which sends the Finsler metric $F_A$ associated with the function $f_A:=f\circ A$ to the Finsler metric $F$. Thus 
\[
\mathrm{sys}(F_A)=\mathrm{sys}(F)=s\quad\text{and}\quad f_A(e_1)=f(u)=s.
\]
The cotangent lift $T^*[A]\colon T^*\T^n\to T^*\T^n$, which has the formula $T^*[A]:=[A]\times (A^{-1})^*$ under the identification $T^*\T^n\cong \T^n\times (\R^n)^*$, is a symplectic automorphism of $(T^*\T^n,\omega_{\mathrm{can}})$ and
\[
T^*[A](D^*_{F_A}\T^n)=D^*_F\T^n.\] In particular, $c(D^*_{F_A}\T^n,\omega_\mathrm{can})=c(D^*_F\T^n,\omega_\mathrm{can})$. Up to relabeling $f_A$ and $F_A$ with $f$ and $F$, this shows that it is enough to prove the upper bound under the assumption that $f(e_1)=s$. 

Every $p\in(\R^n)^*$ has coordinates $p_1:=p\cdot e_1,\ldots,p_n:=p\cdot e_n$. Using these coordinates and reordering the factors, we get a symplectomorphism
\[
\Psi\colon(T^*\T^n,\omega_{\mathrm{can}})\to(\T\times\R,\omega_\mathrm{can})^n,\qquad \Psi(q,p):=\big((q_1,p_1),\ldots,(q_n,p_n)\big).
\]
Since $(f^*)^*=f$, we have for all $k=1,\ldots,n$
\[
s_k:=f(e_k)=\sup_{p\in K^*}|p\cdot e_k|=\sup_{p\in K^*}|p_k|.
\]
Therefore, we obtain
\begin{equation}\label{symemb}
\Psi(D^*_F\T^n)\subset \big(\T\times(-s_1,s_1)\big)\times\ldots\times \big(\T\times(-s_n,s_n)\big).  
\end{equation}
For all $k=1,\ldots,n$, let $r_k:=\sqrt{\frac{2}{\pi}s_k}$ and observe that the maps
\begin{align*}
\phi_k\colon \big (\T \times (-s_k,s_k), \omega_\mathrm{can}\big)\rightarrow (D_{r_k}, \omega_0), \quad \phi_k(q,p)=\sqrt{\frac{p+s_k}{\pi}}e^{2\pi iq}
\end{align*}
are symplectic embeddings \cite[Ch.\ 4.4, Prop.\ 4]{HZ}. By \eqref{symemb}, we get a symplectic embedding
\[
(\phi_1\times\ldots\times  \phi_n)\circ\Psi\colon (D^*_F\T^n,\omega_\mathrm{can})\to (Z_{r_1},\omega_0).
\]
Hence,
\begin{equation*}
c(D^*_F\T^n,\omega_\mathrm{can})\leq \pi r^2_1=2s_1=2s,
\end{equation*}
where the last equality follows from the additional assumption $s_1=f(e_1)=s$.
\subsection{Lower bound}
Denote by $\pi\colon \R^n\to\T^n$ the quotient map. We claim that $\pi$ is injective on $\tfrac{s}{2}K$. Indeed, if $u_1,u_2\in\tfrac{s}{2}K$ and $\pi(u_1)=\pi(u_2)$, then there is $v\in\Z^n$ such that $u_2-u_1=v$. Now $-u_1\in \tfrac s2K$ by reversibility and, thus, $\frac{u_2+(-u_1)}{2}=\tfrac12v\in \tfrac{s}{2}K$ by convexity. Thus, $v\in sK\cap\Z^n$, which implies that $v=0$ by definition \eqref{sys}, and the claim is proved. By the claim, we see that $(\tfrac{s}{2}K)\times K^*\subset D^*_F\T^n$. We get
\begin{equation}\label{fund}
c(D^*_F\T^n,\omega_\mathrm{can})\geq c\big((\tfrac{s}{2}K)\times K^*,\omega_\mathrm{can}\big)=c(K\times K^*,\tfrac{s}{2}\omega_\mathrm{can})=\tfrac{s}{2}c(K\times K^*,\omega_\mathrm{can}). 
\end{equation}
\subsubsection*{Case 1}
Suppose that there is a Finsler norm $G$ induced by a flat Riemannian metric such that $F\geq G$ and $\mathrm{sys}(G)=\mathrm{sys}(F)=s$. In particular, $D^*_G\T^n\subset D^*_F\T^n$ and, therefore, $c(D^*_F\T^n,\omega_\mathrm{can})\geq c(D^*_G\T^n,\omega_\mathrm{can})$. Thus, it is enough to prove the lower bound for $G$, or in other words, we can assume without loss of generality that $F$ is induced by a flat Riemannian metric. In this case, $K$ is an ellipsoid and there is a linear map $A\colon \R^n\to\R^n$ such that $K=A(B_1)$. The map $A\times (A^{-1})^*$ is a symplectomorphism of $(\R^n\times(\R^n)^*,\omega_\mathrm{can})$ and therefore $c(K\times K^*,\omega_\mathrm{can})=c(B_1\times B^*_1,\omega_\mathrm{can})$. By an unpublished argument of Yaron Ostrover and Felix Schlenk \cite[Theorem 12.14(i)]{Sch}, for every $s\in(0,r)$, where $r:=\sqrt{\frac{4}{\pi}}$, there is a symplectic embedding
\[\Phi_s\colon (B_s,\omega_0)\to (B_1\times B_1^*,\omega_\mathrm{can}).
\]
This implies that $c(B_1\times B^*_1,\omega_\mathrm{can})\geq \pi s^2$. Taking the supremum in $s$, we get
\[
c(B_1\times B^*_1,\omega_\mathrm{can})\geq \pi r^2=4.
\]
Using \eqref{fund}, we arrive at the desired lower bound $c(D^*_F\T^n,\omega_\mathrm{can})\geq \tfrac{s}{2}4=2s$.
\begin{Remark}
Although we do not need it here, we observe that there is also a symplectic embedding $\Phi_r\colon (B_r,\omega_0)\to (B_1\times B_1^*,\omega_\mathrm{can})$, see \cite[Lemma 3.7]{Bro} and, for the case $n=2$, \cite{Ram}.
\end{Remark} 
\begin{Remark}
The inequality $c(B_1\times B^*_1,\omega_\mathrm{can})\leq4$ is also true by \cite[Remark 4.2]{AKO}.     
\end{Remark}

\subsubsection*{Case 2}
Suppose that $c\geq c_{\hz}$. We have
\[
c(K\times K^*,\omega_\mathrm{can})\geq c_{\hz}(K\times K^*,\omega_\mathrm{can})=4
\]
by \cite[Theorem 1.7]{AKO}. Using \eqref{fund}, we get $c(D^*_F\T^n,\omega_\mathrm{can})\geq \tfrac{s}{2}4=2s$.
\begin{Remark}
Actually, for every normalized capacity $c$, the inequality $c(K\times K^*,\omega_\mathrm{can})\leq4$ is also true by \cite[Remark 4.2]{AKO}, although we do not need this here. 
\end{Remark}
\bibliographystyle{abbrv}
\bibliography{ref}

\end{document}